\documentclass[a4paper,10pt]{article}
\usepackage[utf8]{inputenc}
\usepackage{color}

\usepackage[all]{xy}
\usepackage{amsmath,amssymb,amsfonts,amsthm,graphicx}
\usepackage[margin=1.2in]{geometry}
\usepackage{multirow}

\newtheorem{thm}{Theorem}[section]

\newtheorem{cor}[thm]{Corollary}
\newtheorem{prop}[thm]{Proposition}

\newtheorem{lm}[thm]{Lemma}

\newtheorem*{thm*}{Theorem}
\newtheorem*{aim*}{Aim}
\newtheorem*{initialaim*}{Initial Aim}
\newtheorem*{conj*}{Conjecture}
\newtheorem*{cor*}{Corollary}
\newtheorem*{prop*}{Proposition}
\newtheorem*{df*}{Definition}
\newtheorem*{lm*}{Lemma}
\newtheorem*{example*}{Example}
\newtheorem*{notation*}{Notation}
\newtheorem*{prob*}{Problem}

\title{On vanishing criteria that control finite group structure}
\author{Julian Brough}

\begin{document}
\date{}
\maketitle

\begin{center}
\small
\textit{Faculty of Mathematics, Centre for Mathematical Sciences,}

\textit{Wilberforce Road, Cambridge, England, CB3 0WA}
\end{center}

\paragraph{}
  \textit{Keywords:}

\textit{Finite groups, vanishing conjugacy classes, soluble groups, super soluble groups}

\normalsize
\begin{abstract}
Many results have been established that show how arithmetic conditions on conjugacy class sizes affect group structure.
A conjugacy class in $G$ is called vanishing if there exists some irreducible character of $G$ which evaluates to zero on the conjugacy class.
The aim of this paper is to show that for some classical results it is enough to consider the same arithmetic conditions on the vanishing conjugacy classes of the group.
\end{abstract}

\section{Introduction}

A well-established area of research in finite group theory considers the relationship between the structure of a group $G$ and sets of positive integers which can naturally be associated to $G$.
One of these sets, the set cs$(G)$ which consists of the conjugacy class sizes of the group $G$, has led to many structural results about $G$.

Some of the classical results concerning cs$(G)$ give arithmetical conditions on cs$(G)$ which yield that $G$ is either a soluble or supersoluble group (see \cite[Theorems 1 and 2]{CWConjCl}).
It is then natural to consider whether the entire data contained in cs$(G)$ is required to obtain such structural results; in fact one of the classical results in this direction shows that if a prime $p$ does not divide the conjugacy class size of any $p'$-element in $G$, then $G$ is $p$-nilpotent (see \cite{AC1}).

More recently, there has been an interest in studying such structural results based on a refinement of the set cs$(G)$ by the character table of the group $G$.
Consider $\chi\in$ Irr$(G)$, an irreducible character of $G$.
A classical result of Burnside says if $\chi$ is non-linear, that is $\chi(1)\not=1$, then there is at least one element $g$ in $G$ such that $\chi(g)=0$ \cite[Theorem 3.15]{IsaChTh}.
In particular, this implies every non-linear row of the character table contains a zero entry.
If one considers conjugacy classes, a natural dual to irreducible characters, then $g$ is a central element in $G$ implies that $|\chi(g)|=\chi(1)$ and thus the column corresponding to $g$ can not contain a zero.
However, it is not clear whether a non-central column must contain a zero.
Therefore we define an element $g$ in $G$ to be a vanishing element if there exists $\chi\in $ Irr$(G)$ such that $\chi(g)=0$.
One can now see that a corollary to Burnside's result is that a group has no vanishing elements if and only if the group is abelian.

Unlike with character values, there is not such a clear statement for a column in the character table to contain a zero.
For nilpotent groups it was shown by Isaacs, Navarro and Wolf that $g$ is non-vanishing if and only if $g$ lies in the centre of $G$ \cite[Theorem B]{INW}. 
For soluble not nilpotent groups they also show if $g$ is non-vanishing then $\overline{g}$ in $G/F(G)$ must be a $2$-element \cite[Theorem D]{INW}. 
In general, if an element $g$ is non-vanishing in $G$ and the order of $g$ is coprime to $6$, then $g$ lies in $F(G)$ \cite[Theorem A]{DNPST}.
Hence one is led to question whether results based on conjugacy class sizes should still hold if one restricts focus only to those corresponding to vanishing elements.

The aim of this paper is to prove three theorems generalising the classical results mentioned above to vanishing conjugacy classes.

\begin{thm*}[Theorem A]
 Let $G$ be a finite group and $p$ a prime dividing the size of $G$ such that if $q$ is any prime dividing the size of $G$, then $q$ does not divide $p-1$.
 Suppose that no vanishing conjugacy class size of $G$ is divisible by $p^2$. Then $G$ is a soluble group.
\end{thm*}

\begin{thm*}[Theorem B]
 Let $G$ be a finite group and suppose that every vanishing conjugacy class size of $G$ is square free.
 Then $G$ is a supersoluble group.
\end{thm*}

\begin{thm*}[Theorem C]
 Let $G$ be a finite group and suppose a prime $p$ does not divide the size of any vanishing conjugacy class size $|x^G|$ for $x$ a $p'$-element of $G$.
 Then $G$ has a normal $p$-complement.
\end{thm*}

In particular, Theorem A has the following corollary.
\begin{cor}
 Let $G$ be a finite group and suppose that no vanishing conjugacy class size of $G$ is divisible by $4$. Then $G$ is a soluble group.
\end{cor}

Theorems A and B form the vanishing equivalents of \cite[Theorem 1]{CWConjCl} and \cite[Theorem 2]{CWConjCl} respectively, while Theorem C is the vanishing analogue of \cite{AC1}.
Note that in the statement of \cite{AC1}, Camina shows that the Sylow $p$-subgroup is a direct factor. 
However, this need not be true when only vanishing elements are considered.
For example the group Sym$(3)$ has a unique vanishing conjugacy class consisting of the transpositions, thus this group has no vanishing $2'$-elements but the Sylow $2$-subgroup is not normal. 
In addition, note the similarity of Theorem C to \cite[Theorem A]{DPSVan}, where the authors showed that if a prime $p$ divides no vanishing conjugacy class size then $G$ has a normal $p$-complement and abelian Sylow $p$-subgroups.

Finally note that these vanishing analogues cover a larger range of groups than the original statements.
Consider the group $C_5\ltimes C_4$, the vanishing elements in this group are the $5'$-elements which each have conjugacy class size equal to $5$.
However the single class of $5$-elements in this group consists of non-vanishing elements but also has conjugacy class size divisible by $4$.
Thus Theorem A implies this group is soluble and Theorem B implies this group is supersoluble, however \cite[Theorem 1]{CWConjCl} and \cite[Theorem 2]{CWConjCl} would not apply to this group as there is a class size divisible by 4.
In addition, Theorem C implies that Sym(3) has a normal $2$-complement, however the conjugacy class of $3$-cycles shows that there exists non-vanishing $2'$-elements with class size divisible by $2$, hence \cite{AC1} can not be applied.

\section{Preliminaries}

Throughout every group considered will be finite.

Recall for $N$ a normal subgroup in $G$, there is a natural bijection between the set of irreducible characters of $G/N$ and the set of irreducible characters of $G$ with $N$ in their kernel.
In particular, this natural bijection implies that if $x$ is an element not in $N$ then $xN$ is vanishing in $G/N$ if and only if $x$ is vanishing in $G$.
In addition, we would like to recall that for an element $x$ in $G$, both $|x^N|$ and $|xN^{G/N}|$ divide $|x^G|$.

\subsection{Vanishing elements in simple groups}

Let $q$ be a prime number, and $\chi$ an irreducible character of $G$; the character $\chi$ is said to have $q$-defect zero if $q$ does not divide $|G|/\chi(1)$.
A result of Brauer highlights the significance $q$-defect zero has for vanishing elements.
If $\chi$ is an irreducible character of $G$ with $q$-defect zero, then $\chi(g)=0$ for every $g\in G$ such that $q$ divides the order of $g$ \cite[Theorem 8.17]{IsaChTh}.

\begin{cor}\cite[Corollary 2]{GODefectZero}\label{ListSimpleCases1}
Let $S$ be a non-abelian simple group and assume there exists a prime $q$ such that $S$ does not have an irreducible character of $q$-defect zero.
Then $q=2$ or $3$ and $S$ is isomorphic either to one of the following sporadic simple groups $M_{12}, M_{22},M_{24},J_2,HS,Suz,Ru,Co_1,Co_3,BM$ or some alternating group $Alt(n)$ with $n\geq 7$.
\end{cor}

In the particular case that $M$ is a minimal normal subgroup, we shall use the preceding corollary with the following lemma.
This result forms a generalisation of a comment made during the proof of \cite[Theorem A]{DPSVan}.

\begin{lm}\label{DefectMinLiftVan}
 Let $G$ be a group, and $N$ a normal subgroup of $G$.
 If $N$ has an irreducible character of $q$-defect zero, then every element of $N$ of order divisible by $q$ is a vanishing element in $G$.
 \begin{proof}
  Take $\psi\in Irr(N)$ of $q$-defect zero.
  Choose $\chi\in Irr(G)$ lying over $\psi$; by Clifford's theorem, Res$^G_N\chi$ is a sum of $G$-conjugates $\psi_i$ of $\psi$.
  As $\psi_i(1)=\psi(1)$, for any $i$, $\psi_i$ is an irreducible character of $N$ with $q$-defect zero.
  Hence each character $\psi_i$ vanishes on any element of $N$ of order divisible by $q$.
  If $x$ is such an element, then $\chi(x)=\sum_i\psi_i(x)=0$, hence $x$ is a vanishing element in $G$.
 \end{proof}
\end{lm}

It still remains to consider those simple groups which have no character of $q$-defect zero for some prime $q$.
The next result provides a condition for an irreducible character of a minimal normal subgroup $M$ of $G$ to extend to an irreducible character of $G$.

\begin{prop}\cite[Lemma 5]{BianChDeg}\label{ExtendingChar}
 Let $G$ be a group, and $M=S_1\times \dots\times S_k$ a minimal normal subgroup of $G$, where every $S_i$ is isomorphic to a non-abelian simple group $S$.
 If $\theta\in Irr(S)$ extends to $Aut(S)$, then $\theta\times\dots\times\theta\in Irr(M)$ extends to $G$.
\end{prop}

The following lemma finds vanishing classes in simple groups which satisfy the condition of Proposition~\ref{ExtendingChar} and also satisfy arithmetical conditions which are required for the proofs of Theorems A and C.

\begin{lm}\label{ListSimpleCases2}
Let $S$ be a non-abelian simple group and assume there exists a prime $q$ such that $S$ does not have an irreducible character of $q$-defect zero:
\begin{enumerate}
\item Then there exists a conjugacy class $x^S$ of size divisible by every prime dividing $S$ and by 4, and there exists $\theta\in Irr(S)$ which extends to $Aut(S)$ such that $\theta$ vanishes on $x^S$.
\item Let $p$ be a prime dividing the order of $S$.
Then there exists a $p'$-element $x$ with conjugacy class $x^S$ of size divisible by $p$,  and there exists $\theta\in Irr(S)$ which extends to $Aut(S)$ such that $\theta$ vanishes on $x^S$.
\end{enumerate}

\begin{proof}
We shall prove $(1)$ and $(2)$ simultaneously.
Observe that if a pair $\{x_1,\theta_1\}$ satisfies the conditions required for $(1)$, then it also satisfies the conditions required for $(2)$, unless $x_1$ turns out to have order divisible by $p$.
Thus to establish $(2)$ from $(1)$, it is enough to provide an additional pair $\{x_2,\theta_2\}$, such that if $x_1$ has order divisible by $p$ then $x_2$ has order not divisible by $p$. 
In particular, the pair $\{x_1,\theta_1\}$ shall be assumed to be the elements chosen in the proof of \cite[Lemma 2.2]{DPSVan}, where the authors provide an element with class size divisible by every prime dividing the size of $S$.

By Corollary~\ref{ListSimpleCases1}, the group $S$ is either isomorphic to a sporadic group $M_{12}$, $M_{22}$, $M_{24}$, $J_2$, $HS$, $Suz$, $Ru$, $Co_1$, $Co_3$, $BM$ or some alternating group Alt$(n)$ with $n\geq 7$.
Assume that $S$ is a sporadic simple group.
For each group the table below provides pairs $\{x_1,\theta_1\}$ and $\{x_2,\theta_2\}$ taken from \cite{Atlas}, as required. 

\begin{center}
\begin{tabular}{lllll}
\hline
Group & Character $\theta_1$ & Class $x_1$ & Character $\theta_2$ & Class $x_2$\\
\hline
$M_{12}$ & $\chi_7$ & $3B$ & $\chi_7$ & $8A$\\
$M_{22}$ & $\chi_3$ & $6A$ & $\chi_2$ & $7A$\\
$M_{24}$ & $\chi_3$ & $6A$ & $\chi_5$ & $7A$\\
$J_2$ & $\chi_6$ & $3B$ & $\chi_{10}$ & $4B$\\
$HS$ & $\chi_7$ & $5C$ & $\chi_{16}$ & $4C$\\
$Suz$ & $\chi_3$ & $8B$ & $\chi_9$ & $3C$\\
$Ru$ & $\chi_2$ & $6A$ & $\chi_9$ & $5B$\\
$Co_1$ & $\chi_2$ & $6H$ & $\chi_2$ & $35A$\\
$Co_3$ & $\chi_9$ & $6E$ & $\chi_{10}$ & $5B$\\
$BM$ & $\chi_2$ & $10D$ & $\chi_6$ & $21A$\\
\end{tabular}
\end{center}

It remains to study the alternating groups.
We provide the conclusion for all $n\geq 7$, although in fact \cite[Corollary 2]{GODefectZero} yields some additional restrictions on $n$.
For $n\geq 7$ recall that Aut(Alt$(n)$)$\cong$Sym$(n)$.

Given a conjugacy class $x^{{\rm Sym}(n)}$ then over Alt$(n)$ this class either stays the same or splits into two equal size conjugacy classes.
Thus for $(1)$ it will be enough to find a conjugacy class in Sym$(n)$ such that $8$ and every odd prime dividing Sym$(n)$ divides $|x_1^{{\rm Sym}(n)}|$. 
While for $(2)$, if $x_1$ has order divisible by $p$, then we need $x_2$ of order not divisible by $p$ such that $p$ or 4 (if $p=2)$ divides $|x_2^{{\rm Sym}(n)}|$.

Every irreducible character of Sym$(n)$ corresponds naturally to a partition of $n$; let $\chi_\sigma$ denote the character corresponding to the partition $\sigma$.
The restriction of $\chi_\sigma$ to Alt$(n)$ is irreducible if and only if the Young diagram corresponding to $\sigma$ is not symmetric \cite{JamKerRepThSym}. 
In the table below pairs $\{x_1,\theta_1\}$ and $\{x_2,\theta_2\}$ are given by a conjugacy class together with a partition $\sigma$, such that $\sigma$ is not symmetric and that it can easily be shown using the Murnaghan-Nakayama formula \cite[Theorem 2.4.7]{JamKerRepThSym} that $\chi_\sigma(x)=0$.

\small
\begin{center}
\hspace*{-1.3cm}
\begin{tabular}{|c|c|c|c|c|c|c|c|}
\hline
\multicolumn{2}{|c|}{$n$} & $\sigma_1$ & Cycle type of $x_1$ & $|x_1^{{\rm Sym}(n)}|$ & $\sigma_2$ & Cycle type of $x_2$ & $|x_2^{{\rm Sym}(n)}|$\\
\hline
\multirow{2}{*}{odd} & not prime & $(n-2,2)$ &  $(n)$ & $(n-1)!$ & \multirow{2}{*}{$(n-3,2,1)$} & \multirow{2}{*}{$(n-2,1,1)$} & \multirow{2}{*}{$n!/2(n-2)$}\\
                                           & prime & $(n-2,2)$ & $(n-3,2,1)$ & $n!/2(n-3)$ & & & \\
\hline
\multirow{2}{*}{even} & $n-1$ not prime & $(n-1,1)$ &  $(n-1,1)$ & $n!/(n-1)$ & \multirow{2}{*}{$(n-2,2)$} & \multirow{2}{*}{$(n-3,1,1,1)$} & \multirow{2}{*}{$n!/6(n-3)$}\\
                                      & $n-1$ prime &  $(n-3,2,1)$ & $(n-2,2)$ & $n!/2(n-2)$ & & & \\
\hline

\end{tabular}
\end{center}

\normalsize

The only case not covered in the above table is when $n=7$ and a conjugacy class size divisible by $4$, because the conjugacy class given above is only divisible by $2$ in Alt$(7)$.
However the character $\chi_2$ and conjugacy class $3B$ taken from \cite{Atlas} satisfy the required conditions of $(1)$ .
\end{proof}
\end{lm}

\subsection{Formations}

The notion of formations provides a powerful tool for studying soluble groups in which every proper quotient has a certain structural property.
First we recall that a set of groups $\mathfrak{F}$ is called a formation if for every $G\in \mathfrak{F}$ then every epimorphic image of $G$ also lies in $\mathfrak{F}$, and if $N_1$ and $N_2$ are normal subgroups of $G$ such that $G/N_1$ and $G/N_2$ both lie in $\mathfrak{F}$, then $G/(N_1\cap N_2)$ lies in $\mathfrak{F}$.
Moreover a formation $\mathfrak{F}$ is called saturated if whenever $G/\Phi(G)$ lies in $\mathfrak{F}$ then $G$ lies in $\mathfrak{F}$, for $\Phi(G)$ the Frattini subgroup of $G$.
In particular, supersoluble groups form a saturated formation \cite[Proposition 10]{HuppertPaper}.

We present the following result, which although standard, we were unable to find a reference for the proof. 
Here $F(G)$ denotes the Fitting subgroup of $G$.

\begin{lm}\label{FGmin}\label{complementMax}
 Let $\mathfrak{F}$ be a saturated formation.
 If $G$ is a soluble group not belonging to $\mathfrak{F}$, but $G/N\in\mathfrak{F}$ for all $N$ normal in $G$, then $F(G)$ is the unique minimal normal subgroup of $G$ and $F(G)$ is an elementary abelian $p$-group for some prime $p$.
Moreover, there exists a complement $H$ to $F(G)$ in $G$ and $H$ is a maximal subgroup of $G$. 
\begin{proof}
 Let $N_1$ and $N_2$ be distinct minimal normal subgroups of $G$, then $G/(N_1\cap N_2)\cong G$ lies in $\mathfrak{F}$.
 Hence $G$ has a unique minimal normal subgroup $N$.

 As $N$ is a $p$-group for $p$ a prime $F(G)$ must also be a $p$-group.
 Furthermore $F(G)/\Phi(F(G))$ is an elementary abelian $p$-group, but as $\Phi(F(G))\leq \Phi(G)=1$, it follows that $F(G)$ is elementary abelian.

Since $N$ is a minimal normal subgroup it has a complement $H$ in $G$ \cite[Hilfssatz VI.7.7]{Huppert}.
 We observe that $N\leq C_G(F(G)\cap H)$ and as $F(G)\cap H\lhd H$, we must have $F(G)\cap H\lhd G$.
 If $F(G)\cap H\not=1$ then $G$ can not have a unique minimal normal subgroup.
It now follows that $N=F(G)$ because $F(G)\cap H=1$.
This proves the first part of the lemma.

Now assume the complement $H$ to $N$ is a proper subgroup of $K$ and $K\leq G$.
As $K\cap N$ is normal in $K$, it follows that $H\leq N_G(K\cap N)$; however $N$ is abelian, so $N\leq N_G(K\cap N)$ and $K\cap N\lhd NH=G$.
Thus $K\cap N=1$ or $N$, but $H<K$ so this is a non-trivial intersection.
Hence $K\cap N=N$ and so $K\leq NH=G$.
\end{proof}
\end{lm}

\section{The proofs}

\begin{thm}[Theorem A]
 Let $G$ be a finite group and $p$ a prime dividing the size of $G$ such that if $q$ is any prime dividing the size of $G$, then $q$ does not divide $p-1$.
 Suppose that no vanishing conjugacy class size of $G$ is divisible by $p^2$. Then $G$ is a soluble group.

 \begin{proof}

Suppose $G$ is chosen of minimal order satisfying the hypothesis of the theorem, but is not soluble.
If $p>2$, then $G$ has odd order and is soluble by the Feit-Thompson Theorem.
Additionally, as the arithmetical conditions on the vanishing classes are inherited by quotients, if $G$ has a proper normal subgroup $N$, then by minimality $G/N$ is soluble.

Let $M$ be a minimal normal subgroup of $G$.
If $M$ is abelian, then $G/M$ is soluble and thus $G$ is soluble.
Hence $M$ must be non-abelian, thus write $M = S_1 \times \cdots \times S_n$, with each $S_i$ isomorphic to a non-abelian simple group $S$.
If $S$ has an irreducible character of $q$-defect zero for every prime $q$, then every non-trivial element of $M$ is a vanishing element in $G$ by Lemma~\ref{DefectMinLiftVan}.
As $S$ is non-soluble, there exists a non-trivial element $x\in S$ such that 4 divides $|x^S|$ by \cite[Theorem 1]{CWConjCl}, and thus $x$ is a vanishing element in $G$ such that 4 divides $|x^G|$.
Now assume there exists a prime $q$ such that $S$ does not have any irreducible characters of $q$-defect zero, then $S$ must be one of the simple groups listed in Corollary~\ref{ListSimpleCases1}.
However, in each case there exists a vanishing element $x$  in the simple group such that 4 divides $|x^S|$ and $x$ vanishes on $\theta\in Irr(S)$ which extends to Aut($S$) by Lemma~\ref{ListSimpleCases2}.
Thus the character $\theta\times\dots\times\theta\in Irr(M)$ extends to $G$ by Proposition~\ref{ExtendingChar}, and $x$ is a vanishing element of $G$ with 4 dividing $|x^G|$.

\end{proof}
\end{thm}

\begin{thm}[Theorem B]
 Let $G$ be a finite group and suppose that every vanishing conjugacy class size of $G$ is square free.
 Then $G$ is a supersoluble group.
\begin{proof}
Suppose $G$ is chosen of minimal order satisfying the hypothesis of the theorem but is not supersoluble; note that $G$ is soluble by Theorem A.
Additionally, as the arithmetical conditions on the vanishing classes are inherited by quotients, if $G$ has a proper normal subgroup $N$ then by minimality $G/N$ must be supersoluble.

As the set of all supersoluble groups form a saturated formation \cite[Proposition 10]{HuppertPaper}, $G$ must have a unique minimal normal subgroup $N=F(G)$ by Lemma~\ref{FGmin}. 
Moreover $N$ must have order $p^k>p$, as if $N$ had order $p$ then $G$ would be supersoluble by induction on $G/N$.
In addition, there exists a complement $H$ in $G$ for $N$ \cite[Hilfssatz VI.7.7]{Huppert} and $G/N\cong H$ must be supersoluble.

Let $M$ be a minimal normal subgroup of $H$.
Then $|M|=q$ for some prime $q$, and as $O_p(G)=F(G)$, it follows that $q\ne p$.
For any $x\in M\backslash \{1\}$, the subgroup $C_G(x)$ is contained in $N_G(M)=N_G(\langle x\rangle)$.
As $H\leq N_G(M)$, it follows that $N_G(M)=H$ or $G$ by Lemma~\ref{complementMax}.
Thus $N_G(M)=H$ and hence $|N|$, and therefore $p^2$, divides the size of $|x^G|$.
If $x$ is a vanishing element in $G$, then this yields a contradiction.

Thus assume for any minimal normal subgroup $M$ of $H$, the elements $x\in M\backslash \{1\}$ are non-vanishing in $G$.
As $H$ naturally maps isomorphically to $G/F(G)$, by sending $x$ in $H$ to $\overline{x}=xF(G)$, the order of $x$ in $H$ equals the order of the image of $x$ in $G/F(G)$, which by \cite[Theorem D]{INW} is a power of $2$.
In particular, as $M=\langle x\rangle$ and the size of $M$ is prime, every minimal normal subgroup of $H$ must be isomorphic to $C_2$.
If the largest prime $q$ which divides $|H|$ is odd, then $H$ must have a minimal normal subgroup of order $q$ \cite[Lemma 19.3.1]{Hall}, contradicting that $M\cong C_2$.
Hence $H$ must be a $2$-group, and is therefore nilpotent.
Therefore $F(G)$ is the penultimate term of the ascending fitting series of $G$ and all non-vanishing elements of $G$ lie in $ F(G)=N$ \cite[Theorem 2.4]{INW}, which implies $x\in M\backslash 1$ is vanishing in $G$.
Thus $x$ must be vanishing in $G$, yielding a contradiction.

\end{proof}
\end{thm}

\begin{thm}[Theorem C]
 Let $G$ be a finite group and suppose a prime $p$ does not divide the size of any vanishing conjugacy class size $|x^G|$ for $x$ a $p'$-element of $G$.
 Then $G$ has a normal $p$-complement.
 \begin{proof}
  Suppose $G$ is chosen of minimal order satisfying the hypothesis of the theorem, but does not have a normal $p$-complement.
Additionally, as the arithmetical conditions on the vanishing classes are inherited by quotients, if $K=O_{p'}(G)\not=1$, then by minimality $G/K$ has a normal $p$-complement and hence so does $G$.
  Thus assume that $ O_{p'}(G)=1$.

 Let $M=S_1\times\dots\times S_k$ be a minimal normal subgroup of $G$, with each $S_i\cong S$ a simple group.
As $ O_{p'}(G)=1$ it follows that $p$ divides the order of $M$.
 Assume first that $S$ is non-abelian.
 If $S$ has an irreducible character of $q$-defect zero for every prime $q$, then by Lemma~\ref{DefectMinLiftVan} every non-trivial element of $M$ is a vanishing element in $G$.
As $S$ has no normal $p$-complement, by \cite{AC1}, there exists a $p'$-element $x$ in $S$ such that $p$ divides $|x^S|$, and thus $x$ is a $p'$-element which is vanishing in $G$ such that $p$ divides $|x^G|$.

 Now assume there exists a prime $q$ such that $S$ does not have any irreducible characters of $q$-defect zero, then $S$ must be one of the simple groups in Corollary~\ref{ListSimpleCases1}.
 However, by Lemma~\ref{ListSimpleCases2} there exists a $p'$-element $x$ such that $p$ divides $|x^S|$ with a corresponding character $\theta\times\dots\times\theta\in Irr(M)$ which extends to $G$ by Proposition~\ref{ExtendingChar}.
 Therefore $x$ is a vanishing $p'$-element of $G$ with $p$ dividing $|x^G|$.

 Finally consider the case that $S$ is an abelian $p$-group and so $N:=O_p(G)\not=1$.
 By the minimality of $G$, the quotient $G/N$ has a normal $p$-complement so $G/N$ is $p$-soluble; hence $G$ is also $p$-soluble.
 As $O_{p'}(G)=1$, then by \cite[Theorem 6.3.2]{Gor} $N$ is a self centralising subgroup of $G$, i.e. $C_G(N)\leq N$.
 Thus for $g\in G$ if $p$ does not divide $|g^G|$ then $g\in N$ and hence any vanishing $p'$-element of $G$ lies in $N$.
 In particular, there are no vanishing elements of $p'$ order in $G$.
 Hence $G$ has a normal $p$-complement by \cite[Corollary C]{DPSVanOrd}.

 \end{proof}

\end{thm}

The author would like to thank Emanuele Pacifici and Mariagrazia Bianchi for their hospitality during his stay at the Universit\`{a} Degli Studi Di Milano in September 2014, where he was first introduced to the concept of vanishing conjugacy classes.

\bibliographystyle{plain}
\bibliography{bibfile}

\end{document}